\magnification=1070
\overfullrule=0mm

\input xy \xyoption{all}

\font\tenrom=cmr10 
\font\fiverom=cmr10 at 5pt
\font\sevenrom=cmr10 at 7pt

\font\nineromi=cmmi10 at 9pt

\font\sevencmb=cmb10 at 7pt

\font\ninerom=cmr10 at 9pt
\font\nineromtt=cmtt10 at 9pt
\newfam\romfam
\scriptscriptfont\romfam=\fiverom
\textfont\romfam=\tenrom
\scriptfont\romfam=\sevenrom
\def\rom{\fam\romfam\tenrom}

\font\oldstylen=cmmi10 at 9pt

\def\ZBIc{\hbox{\sevencmb Z}}
 
\newtoks\auteurcourant      \auteurcourant={\hfil}
\newtoks\titrecourant       \titrecourant={\hfil}

\newtoks\hautpagetitre      \hautpagetitre={\hfil}
\newtoks\baspagetitre       \baspagetitre={\hfil}

\newtoks\hautpagegauche   \newtoks\hautpagedroite 
  
\hautpagegauche={\eightcmr\rlap{\folio}\eightcmr\hfil\the\auteurcourant\hfil}
\hautpagedroite={\eightcmr\hfil\the\titrecourant\hfil\eightcmr\llap{\folio}}

\newtoks\baspagegauche      \baspagegauche={\hfil} 
\newtoks\baspagedroite      \baspagedroite={\hfil}

\newif\ifpagetitre          \pagetitretrue

\headline={\ifpagetitre\the\hautpagetitre
            \else\ifodd\pageno\the\hautpagedroite
             \else\the\hautpagegauche
              \fi\fi}

\footline={\ifpagetitre\the\baspagetitre\else
            \ifodd\pageno\the\baspagedroite
             \else\the\baspagegauche
              \fi\fi
               \global\pagetitrefalse}

\def\raggedbottom{\topskip 10pt plus 36pt\r@ggedbottomtrue}

\newcount\notenumber \notenumber=1
\def\note#1{\footnote{$^{{\the\notenumber}}$}{\eightcmr {#1}}
\global\advance\notenumber by 1}

\auteurcourant={Philippe NUSS -- Marc WAMBST}
\titrecourant={NON-ABELIAN HOPF COHOMOLOGY}

\font\bb=msym10

\font\tengoth=yfrak at 12pt
\font\fivegoth=yfrak at 5pt
\font\sevengoth=yfrak at 7pt
\newfam\gothfam
\scriptscriptfont\gothfam=\fivegoth
\textfont\gothfam=\tengoth
\scriptfont\gothfam=\sevengoth
\def\goth{\fam\gothfam\tengoth}
\font\eightcmr=cmr10 at 8pt

\font\ninecmsl=cmsl10 at 9pt
\font\ninecmbx=cmbx10 at 9pt
\font\twelvecmbx=cmbx10 at 12pt

\font\tensymb=msxm9
\font\fivesymb=msxm9 at 5pt
\font\sevensymb=msxm9 at 7pt
\newfam\symbfam
\scriptscriptfont\symbfam=\fivesymb
\textfont\symbfam=\tensymb
\scriptfont\symbfam=\sevensymb

\font\tensymbo=msym9
\font\fivesymbo=msym9 at 5pt
\font\sevensymbo=msym9 at 7pt
\newfam\symbofam
\scriptscriptfont\symbofam=\fivesymbo
\textfont\symbofam=\tensymbo
\scriptfont\symbofam=\sevensymbo

\font\tensymbol=cmsy10
\font\fivesymbol=cmsy10 at 5pt
\font\sevensymbol=cmsy10 at 7pt
\newfam\symbolfam
\scriptscriptfont\symbolfam=\fivesymbol
\textfont\symbolfam=\tensymbol
\scriptfont\symbolfam=\sevensymbol

\newfam\bffam
\scriptscriptfont\bffam=\fivebf
\textfont\bffam=\tenbf
\scriptfont\bffam=\sevenbf
\def\bf{\fam\bffam\tenbf}

\newfam\mathb

\def\Z{\hbox{\bb Z}}
\def\Zb{\hbox{\bf Z}}

\def\Z2{\hbox{$\Zb/2\Zb$}}

\def\Cri{{\mathop{\rm C}\nolimits}}
\def\D{\hbox{$\cal D$}}

\def\T{\hbox{$\cal T$}}
\def\U{\hbox{$\cal U$}}
\def\V{\hbox{$\cal V$}}

\def\Cr{\hbox{\rm C}}
\def\Dr{\hbox{\rm D}}
\def\Hr{\hbox{\rm H}}
\def\Wr{\hbox{\rm W}}
\def\Zr{\hbox{\rm Z}}

\def\adj{\hbox{\tensymbol{\char 97}}}

\def\Del {\hbox{$\Delta$}}

\def\id{\hbox{\rm id}}

\def\eps {\hbox{$\varepsilon$}}

\def\dm{ \null\hfill {\fivesymb {\char 3}}}
\def\lr {\hbox{$\ \longrightarrow\ $}}

\def\ot {\hbox{$\otimes$}}
\def\pa{\S\kern.15em}
\def\Dem{\noindent {\sl Proof.}$\, \,$}

\def\longrighthook{\lhook\joinrel\relbar\joinrel\rightarrow}

\def\square{\hbox{ {$\sqcap $} \kern -1.315em
\raise -0,3325ex \hbox { $\scriptstyle -$}$\!$ }}

\def\builda#1_#2^#3{\mathrel{\mathop{\kern 0pt#1}\limits_
{\kern -2.26em{#2}}^{\kern -2.26em{#3}}}}

\def\hflong{\smash{\mathop{\hbox to 10mm{\rightarrowfill}}
\limits}}

\def\dblarrow{\raise 0,4ex \hbox{ {$\longrightarrow $} \kern -2.26em
\raise -0,7ex \hbox { $\longrightarrow$}$\!$ }}

\def\dblbiarrow{\raise 0,4ex \hbox{ {$\longrightarrow $} \kern -2.26em
\raise -0,7ex \hbox { $\longleftarrow$}$\!$ }}

\def\builda#1_#2^#3{\mathrel{\mathop{\kern 0pt#1}\limits_
{\kern -2.26em{#2}}^{\kern -2.26em{#3}}}}

\def\hfl#1#2{\smash{\mathop{\hbox to 8mm{\rightarrowfill}}
\limits^{\scriptstyle#1}_{\scriptstyle#2}}}

\def\hfll#1#2{\smash{\mathop{\hbox to 10mm{\rightarrowfill}}
\limits^{\scriptstyle#1}_{\scriptstyle#2}}}

\def\hflll#1#2{\smash{\mathop{\hbox to 25mm{\rightarrowfill}}
\limits^{\scriptstyle#1}_{\scriptstyle#2}}}

\def\lhfll#1#2{\smash{\mathop{\hbox to 10mm{\leftarrowfill}}
\limits^{\scriptstyle#1}_{\scriptstyle#2}}}

\def\trplarrow{\raise 0,4ex \hbox{ {$\dblarrow$} \kern -2.45em
\raise -1,0ex \hbox { $\longrightarrow$}$\!$ }}

\def\qplarrow{\raise 0,4ex \hbox{ {$\dblarrow $} \kern -2.8em
\raise -1,42ex \hbox { $\dblarrow$}$\!$ }}

\def\Aut{{\rm Aut}}
\def\End{\mathop{\rm End}\nolimits}
\def\Tors{{\rm Tors}}
\def\tors{{\rm tors}}

\def\relbar{\mathrel{\smash-}}

\def\lodot{ \mathop{\circ\mkern-7.05mu\cdot\mkern3.5mu}\nolimits }

\noindent {\twelvecmbx NON-ABELIAN HOPF COHOMOLOGY

}
\smallskip
\noindent {\twelvecmbx 
}
\vskip 25pt

\noindent {PHILIPPE NUSS, MARC WAMBST}
\vskip 3pt
\noindent {\ninerom Institut de Recherche Math\'ematique Avanc\'ee,
Universit\'e Louis-Pasteur et CNRS, 7, rue Ren\'e-Descartes,
67084 Strasbourg Cedex, France. 
e-mail: {\nineromtt nuss@math.u-strasbg.fr} and {\nineromtt wambst@math.u-strasbg.fr}}

\vskip 20pt
\noindent {\ninecmbx Abstract.} {\ninerom We introduce non-abelian cohomology sets of Hopf algebras with coefficients in Hopf modules.
We~prove that these sets generalize Serre's non-abelian group cohomology theory. Using descent techniques, we esta\-blish that 
our construction enables to classify as well twisted 
forms for modules over Hopf-Galois extensions as torsors over  Hopf-modules.}


\vskip 5pt
\noindent {\ninecmbx Key-words:} {\ninerom non-abelian cohomology, noncommutative descent theory, Hopf-Galois extension,
Hopf-module,
twis\-ted form, torsor, Hilbert's Theorem~90}

\vskip 30pt
\noindent {{\bf I{\ninecmbx NTRODUCTION}.} 
The aim of this article is to extend to Hopf algebras the concept of 
non-abelian cohomology of groups. Introduced in 1958 by Lang and Tate ([8])
for Galois groups with coefficients in an algebraic group, the non-abelian cohomology 
theory in degree 0 and 1 was formalized by Serre ([12], [13]).
For an arbitrary group $G$ acting
on a  group  $A$ which is not necessarily abelian, Serre constructs 
a $0$-cohomology group $\Hr^0(G, A)$ and a $1$-cohomology pointed set $\Hr^1(G, A)$.
These objects generalize the
two first groups of the classical  Eilenberg-MacLane
cohomology sequence  $\Hr^*(G, A) = {\rm Ext}^*_{{\hbox {\ZBIc}} [G]}(\Zb, A)$, defined only when 
 $A$ is abelian. It is well-known that the non-abelian cohomology set $\Hr^1(G, A)$ classifies
the torsors on $A$ (see [13]).

The non-abelian cohomology theory of groups comes naturally into play in the particular case
where $S/R$ is a $G$-Galois extension of rings in the sense of [9]. The situation is the following:
a finite group $G$
acts on a ring extension $S/R$ and, in a compatible way, on an $S$-module $M$. The coefficient group is then
the group of $S$-automorphisms $A = \Aut _S(M)$ of  $M$. In [10],
one of the authors showed  that the set $\Hr^1(G, \Aut _S(M))$ classifies as well descent cocycles 
on $M$ as twisted forms of $M$. 

Galois extensions of rings may be viewed as particular cases
of Hopf-Galois extensions defined by Kreimer-Takeuchi ([7]), where a Hopf algebra $H$ 
(non necessarily commutative nor cocommutative) plays the
r\^ole of the Galois group. Indeed, given a group  $G$, a $G$-Galois extension of rings is nothing but a $\Zb^G$-Hopf-Galois
extension of rings, where $\Zb^G$ stands for the dual Hopf algebra of the group ring $\Zb[G]$.

\medskip
Suppose now  fixed a ground ring $k$, a Hopf algebra $H$ over  $k$, and an
$H$-comodule algebra $S$ (for instance, any $H$-Hopf-Galois extension $S/R$ is based on such a datum).
For any $(H,S)$-Hopf module $M$,
that is an abelian group $M$ endowed with an  $S$-action  and a compatible $H$-coaction,
we define in the cosimplicial spirit a $0$-cohomology group $\Hr^0(H, M)$ and
a $1$-cohomology pointed set  $\Hr^1(H, M)$.

The philosophy behind the construction is the following (precise definitions will be given in the core of the paper). 
Start with a $G$-Galois extension $S/R$, where $G$ is a finite group, and with $M$  a  $(G,S)$-Galois module,
{\sl i.e.} an abelian group $M$ endowed with two compatible $S$- and $G$-actions. The group $\Aut _S(M)$ 
inherits a $G$-action by conjugation.  
Let $k^G$ be the    dual Hopf algebra of the group ring $k[G]$.
A $1$-cocycle in the sense
of Serre is represented by a certain map $\alpha: G \lr \Aut _S(M)$.  
By duality, $\alpha$ formally defines
an element in $M \ot_k M^* \ot_k k^G$, which can also be seen as a map $\Phi _\alpha: M \lr M \ot _k k^G$
satisfying some conditions.
Assume now given, instead of $G$, a Hopf-algebra $H$ coacting on a ring~$S$.
Let $M$ be an $(H,S)$-Hopf module, that is a module on which both $H$ and $S$ act in a compatible way. 
We replace the former map $\Phi _\alpha: M \lr M \ot _k k^G$ by a map $\Phi: M \lr M \ot _k H$
and state general
requirements -- the cocycle conditions --, which  reflect the group-cocycle condition
on $\alpha$. 
This construction gives rise to a $1$-cohomology pointed set $\Hr^1(H, M)$.

\medskip
We establish two mains results. The first Theorem shows that the $1$-cohomology set $\Hr^1(H, M)$
generalizes the non-abelian group $1$-cohomology set of Serre. The second one relates $\Hr^1(H, M)$
to ${\rm Twist }(S/R, N_0)$, the isomorphy class of the twisted forms of an extended module  $M = N_0\ot_RS$.
More precisely, we prove the two following statements:

\medskip
\noindent {\bf Theorem A.} {\sl  For  a group $G$ and a $(k^G,S)$-Hopf module $M$, there is an isomorphism of pointed sets
$$ \Hr^1(k^G, M) \cong \Hr^1(G, \Aut _S(M)).$$ }

\noindent {\bf Theorem B.} {\sl  For a Hopf-algebra $H$  and an $(H,S)$-Hopf module  $M$  of the form $M = N_0\ot_RS$, 
there is an isomorphism of pointed sets
$$ \Hr^1(H, M) \cong {\rm Twist }(S/R, N_0).$$}

\noindent  The precise wording of Theorem A
will be found in Theorem 3.2, and that of Theorem B in Theorem~1.2.
As a consequence of Theorem B, we deduce (Corollary 1.3) a Hopf version of the celebrated Theorem 90
stated in 1897 by Hilbert in his {\sl  Zahlbericht}.

\medskip

In order to prove these two results, we bring in an auxiliary cohomology  theory $\Dr^i(H, M)$ ($i=0,1$)
related to Descent Theory. The pointed set  $\Dr^1(H, M)$ classifies the $(H,S)$-Hopf
{\hbox{module}} structures on $M$ and, in the case of a Hopf-Galois extension, the descent data on $M$.
Moreover, it may be viewed as torsors on $M$ (Proposition 2.8).
%It is obtained by taking the subset of  $\Hom _k(M, M\ot _k H)$ consisting of Hopf-Galois descent data 

\medskip
We mention here that A. Blanco Ferro ([1]), generalizing a construction due to M. Sweedler
([14]), defined a
$1$-cohomology set $\Hr^1(H, A)$, where $H$ is a Hopf-algebra and $A$ is an $H$-module algebra. He~applied his
theory, which is in some sense dual to ours, to a commutative particular case: 
not only does $H$ have to be a commutative finitely generated $k$-projective Hopf algebra,
but $S/k$ is a commutative Hopf-Galois extension. For any $k$-module $N$,
setting $A = \End _S(N\ot _kS)$, Blanco Ferro showed in this particular case that his set
$\Hr^1(H^*, A)$ classifies the twisted forms of $N\ot _kS$ where $H^*$ stands for the dual Hopf algebra of $H$.

\smallskip
\noindent

\medskip
\medskip
\noindent {\bf 0. Conventions.}

\smallskip

\noindent 
Let $k$ be a fixed commutative and unital ring.
%%which will in some results be a field. 
The unadorned symbol $\ot$ between a right $k$-module and a left
$k$-module stands for $\ot _k$.
By {\sl  algebra} we mean a unital associative $k$-algebra.
A {\sl  division algebra} is either a commutative field or a skew-field. 
By {\sl  module} over a ring $R$, we always understand a right $R$-module unless otherwise stated.
Denote by ${\goth Mod}_R$
the  category of $R$-modules and by ${\goth Set}$ the category of sets.

Let $H$ be a finite-dimensional Hopf-algebra over $k$ with multiplication $\mu_H$, unity map $\eta _H$,
comultiplication $\Delta _H$, counity map $\varepsilon _H$, and  antipode $\sigma_H$.
%%We assume that the antipode 
%%$\sigma_H$ of $H$ is invertible.
Let $S$ be an algebra, $\mu _S$ its multiplication, $\eta _S$ its unity map. We assume that 
$S$ is  a right $H$-comodule algebra, in other words that $S$ is equipped with an $H$-coaction map $\Delta _S: S \lr S \ot H$
which is a morphism of algebras.
Let $M$ be both an $S$-module 
and an $H$-comodule with the $H$-coaction map
$\Delta _M: M \lr M \ot H$. If $\Delta _M$ verifies the equality
$$\Delta _M (ms) = \Delta _M(m) \Delta _S(s), \eqno (1)$$ for any $m \in M$ and $s \in S$, we say that
$M$ is an {\sl  $(H,S)$-Hopf module} 
(also called a {\sl  relative Hopf module} in the literature)  and that
$\Delta _M: M \lr M \ot H$ is  {\sl  $(H,S)$-linear}.
A morphism $f: M \lr M'$~of {\hbox{$(H,S)$-Hopf}} modules is an $S$-linear map
$f$ such that $(f \ot \id _M) \circ\Delta _M = \Delta _{M'} \circ f$.
%%Denote by ${\goth Mod}_S^H$ the category of 
%%$(H , S)$-Hopf modules.
To denote the coactions on elements, we use the Sweedler-Heyneman convention, that is, for $m \in M$, we write
 $\Delta _M(m) = m_0 \ot m_1$, with summation implicitly
understood. More generally, when we write down a tensor we usually omit the summation sign $\sum$.

\medskip
Denote by $R$ the algebra of $H$-coinvariants of $S$, that is $R = \{ s \in S \ \vert \ \Delta _S (s) = s\ot 1\}$.
An $S$-module $M$ is said to be {\sl  extended} if there exists an $R$-module $N_0$ such that
$M$ is equal to $N_0 \ot _RS$.
The inclusion map $\psi: R \longrighthook S$ 
is a {\sl  (right) $H$-Hopf-Galois
extension} if $\psi$ is faithfully flat 
and the map $\Gamma _{\psi}: S \ot _RS \lr S \ot H $, called {\sl  Galois map}, given on an indecomposable tensor $s \ot t \in S \ot _RS$ by
$$\Gamma _{\psi}(s\ot t) = s\Del _S(t),$$ is a $k$-linear isomorphism. 
By Hopf-Galois descent theory ([5], [11]), every $(H,S)$-Hopf module is isomorphic to an extended $S$-module.
Conversely, an extended $S$-module $M = N_0 \ot _RS$ owns an  $(H,S)$-Hopf module structure 
with the canonical coaction $\Delta _M  = \id _{N_0} \ot \Delta _S: N_0 \ot _RS \lr N_0 \ot _RS \ot H$.
%%Call $\Gamma_{\psi}$  the {\sl  Hopf-Galois isomorphism of $\psi$}.
\medskip

Let $G$ be a finite group. Denote by $k^G$ the $k$-free Hopf algebra over the $k$-basis
$\{ \delta _g\} _{g \in G}$, with the following structure maps:
the multiplication is given by ${\displaystyle \delta _g \cdot \delta _{g'} = \partial _{g,g'} \delta _{g}}$,
where $\partial _{g,g'}$ stands for the 
Kronecker symbol of $g$ and $g'$; 
the  comultiplication $\Del _{k^G}$ is
defined by 
${\displaystyle \Del _{k^G} (\delta _g) =
\sum _{ab = g}\delta _a \ot \delta _b}$;
the unit in $k^G$ is the element 
${\displaystyle 1 = \sum _{g \in G}\delta _g}$;
the counit $\varepsilon _{k^G}$ is
defined by $\varepsilon _{k^G} (\delta _g) = \partial _{g,e}1$; 
the antipode $\sigma_{k^G}$ sends $\delta _g$ on $\delta _{g^{-1}}$.
When $k$ is a field, then $k^G$ is the dual of the usual group Hopf-algebra $k[G]$.
It is easy to see that a $k^G$-Hopf-Galois extension is the same as a $G$-Galois extension of $k$-algebras
in the sense of [9].
To give an action of $G$ on $S$ is equivalent to give a coaction map of $k^G$ on $S$, the two structures being related by
the equality

$$\Delta _S(s) = \sum _{g \in G} g(s) \ot \delta _g.$$

An $S$-module $M$ will be called a {\sl  $(G,S)$-Galois module} if it is endowed with a $(G,S)$-action, that is a
$G$-action $\gamma: G  \lr \Aut _k(M)$ such that following twisted $S$-linearity condition:
$$g(ms) = g(m)g(s) \eqno (2)$$ holds for any $g \in G$, $m\in M$, and $s \in S$
(when no confusion about $\gamma$ is possible, we denote for simplicity $g(m)$ instead of $\gamma (g)(m)$). 
When $\gamma$ verifies (2), we say that the morphism $\gamma$ is {\hbox{{\sl  $(G,S)$-linear}}}.
Denote by $\Aut^\gamma _{S}(M)$ the subgroup  of $\Aut _k (M)$ which is the image of $\gamma$.

To give a $(G,S)$-Galois module structure on $M$ is equivalent to give a $(k^G,S)$-Hopf
module structure on $S$. By Galois descent theory, a  $(G,S)$-Galois module is isomorphic to an
extended module $N \ot _RS$.

%%generated
%%by the $(G,S)$-actions. 

\bigskip

\noindent {\bf 1. Non-abelian Hopf cohomology theory.}

\smallskip
\noindent In this section we define a non-abelian Hopf cohomology theory,  and state our main result, Theorem~1.2, which
compares in the Hopf-Galois context the $1$-Hopf cohomology set 
with twisted forms. We deduce a Hopf-Galois version
of Hilbert's Theorem 90.

\smallskip

\noindent {\sl  1.1. Definition of the non-abelian Hopf cohomology sets.}

\smallskip

\noindent Let $H$ be a Hopf-algebra and $S$ be an $H$-comodule algebra. For any $S$-module
 $M$, we endow $M \ot H^{\otimes n}$ with an $S$-module structure given by
$$(m \ot \underline {h}) s = ms \ot \underline {h},$$ for $m \in M$, $\underline {h} \in H^{\otimes n}$, and
$s \in S$. 
\goodbreak
Set
$\Wr_k^n(M) = {\rom Hom}_k(M , M \ot H^{\otimes n})$ and
$\Wr_S^n(M) = {\rom Hom}_S(M , M \ot H^{\otimes n})$. We equip the {\hbox {$k$-module}} $\Wr_k^n(M)$ 
with a composition-type product $\lodot: \Wr_k^n(M) \ot \Wr_k^n(M) \lr \Wr_k^n(M)$,
defined by
$$\left\{\eqalign{\varphi \lodot \varphi' & = \varphi \circ \varphi' \ \ \ {\rom if} \ \ n = 0\cr
\varphi \lodot \varphi' & = (\id _M \ot \mu _H^{\otimes n})\circ (\id _M \ot \chi_n) \circ (\varphi \ot \id _H^{\otimes n}) 
\circ \varphi' \ \ {\rom if} \ \ n > 0\cr}\right.$$ 
for $\varphi, \varphi' \in \Wr_k^n(M)$;
here $\chi_n: H^{\otimes n}\ot H^{\otimes n} \lr (H\ot H)^{\otimes n}$ denotes the intertwining operator given by
$$\chi_n\bigl((a_1 \ot \ldots \ot a_n) \ot (b_1 \ot \ldots \ot b_n)\bigr) = (a_1 \ot b_1) \ot \ldots \ot (a_n \ot b_n).$$
It restricts to a product still denoted $\lodot$ on $\Wr_S^n(M)$.
Thanks to the product $\lodot \hskip-1pt$, the modules $\Wr_k^n(M)$ and $\Wr_S^n(M)$ 
become a monoid: the associativity of $\lodot$ is a direct
consequence of the coassociativity of $\Delta _H$ and the neutral element is $\upsilon_n = \id _M \ot \eta _H^{\otimes n}$.
Further we shall use that   the group of invertible elements of the monoid  $\Wr_S^0(M)$ is $\Aut _S(M)$.

\medskip

Suppose that $M$ is an $H$-comodule. Denote by $T$ the {\sl  flip} of $H \ot H$, the automorphism of $H \ot H$ which sends an indecomposable tensor $h \ot h'$ to
$h'\ot h$.
We define two maps $d^i: \Wr_k^0(M) \lr \Wr_k^1(M)$
($i= 0,1$)
and three maps $d^i: \Wr_k^1(M) \lr \Wr_k^2(M)$ ($i=0,1,2$) %and $\delta: M \lr M \ot H \ot H$
by the formulae 
$$\eqalign{ d^0\varphi   & = (\id_{M} \ot \mu _H) \circ (\Delta _{M} \ot \id _H)
\circ (\varphi \ot \sigma _H) \circ \Delta _M \cr 
d^1\varphi & = (\id _{M} \ot \eta _H)\circ \varphi\cr
 d^0\Phi & = (\id _{M} \ot \mu _H \ot \id _H)
 \circ (\Delta _{M} \ot T) \circ (\Phi  \ot \sigma _H) \circ \Delta _{M}\cr
d^1\Phi & = (\id _{M} \ot  \Delta _H) \circ \Phi \cr
d^2\Phi & = (\id _{M} \ot \id _H \ot \eta _H) \circ \Phi  = \Phi \ot \eta _H, \cr
}$$
where $\varphi: M \lr M$ and $\Phi: M \lr M \ot H$ are  $k$-linear morphisms.
\medskip

\noindent {\bf Lemma 1.1.} {\sl  Let $M$ be an $(H,S)$-Hopf-module.
The restriction of the above defined maps to the corresponding
monoids $ \Wr_S^0(M)$ and $ \Wr_S^1(M)$ are morphims of monoids which may be organized in the following cosimplicial 
diagram:

$$\xymatrix{ \Wr_S^0(M)  \ar@<1.3ex>[r]^{d^0} 
\ar@<-1.3ex>[r]^{d^1}  &  \Wr_S^1(M) 
\ar@<2.3ex>[r]^{d^0} 
\ar@<0ex>[r]^{d^1} 
\ar@<-2.3ex>[r]^{d^2}&  \Wr_S^2(M)} \eqno (3)$$}

\Dem We adopt  the Sweedler-Heyneman convention and use
the Hopf yoga, for instance, the fact that for any $x,y \in H$,
one has $x_0\ot \sigma_H(x_1)x_2y = x_0\ot \varepsilon_H(x_1)y = x \ot y$. 
 First one has to show that $d^i\varphi$ and $d^i\Phi$ are $S$-linear. This assertion is obvious for 
$d^1\varphi$. Let us prove it for  $d^0\varphi$. We get,
for any $m\in M$ and $s \in S$, the equalities
$$\eqalign { d^0\varphi (ms) & = 
[(\id_{M} \ot \mu _H) \circ (\Delta _{M} \ot \id _H) \circ (\varphi \ot \sigma _H) \circ \Delta _M](ms) \cr
&=
[(\id_{M} \ot \mu _H) \circ (\Delta _{M} \ot \id _H)](\varphi(m_0)s_0 \ot \sigma _H(m_1s_1)) \cr
&= (\id_{M} \ot \mu _H) [(\varphi(m_0)_0s_0 \ot \varphi(m_0)_1s_1 \ot\sigma _H(s_2)\sigma _H(m_1)] \cr
& = \varphi(m_0)_0s_0 \ot \varphi(m_0)_1(s_1\sigma _H(s_2))\sigma _H(m_1)  \cr
& = \varphi(m_0)_0s \ot \varphi(m_0)_1\sigma _H(m_1) \cr & = d^0\varphi (m)s.\cr}$$
The $S$-linearity of $d^1\Phi$ and $d^2\Phi$ is obvious. 
We prove it for  $d^0\Phi$. For any $m\in M$ and $s \in S$, set
$\Phi(m) = m' \ot m''$. We  have $d^0\Phi (m) = ((m_0)')_0\ot ((m_0)')_1\sigma _H(m_1)\ot (m_0)''$, hence
$$\eqalign { d^0\Phi (ms) & =  [(\id _{M} \ot \mu _H \ot \id _H)
 \circ (\Delta _{M} \ot T) \circ (\Phi  \ot \sigma _H) \circ \Delta _{M}](ms) \cr
&=
 [(\id _{M} \ot \mu _H \ot \id _H)
 \circ (\Delta _{M} \ot T)] \bigl((m_0)'s_0\ot (m_0)'' \ot \sigma _H(m_1s_1)\bigr)\cr
&= (\id _{M} \ot \mu _H \ot \id _H)
 [((m_0)')_0s_0\ot ((m_0)')_1s_1 \ot \sigma _H(s_2)\sigma _H(m_1)\ot (m_0)''] \cr
& = ((m_0)')_0s\ot ((m_0)')_1\sigma _H(m_1)\ot (m_0)'' \cr
& = d^0\Phi (m) s.\cr}$$

We prove now that $d^i$ respects the monoid structures on $\Wr_S^k(M)$,
that is  $$d^i\varphi \lodot d^i\varphi ' = d^i(\varphi \lodot \varphi ') \hbox{\rom,}\ \ \
d^i\Phi \lodot d^i\Phi ' = d^i(\Phi \lodot \Phi '),
\ \ \ \hbox{\rom and}\ \ \ d^i(\upsilon_k) = \upsilon_{k+1}$$
for any $\varphi , \varphi' \in \Wr_S^0(M)$, any $\Phi , \Phi' \in \Wr_S^1(M)$, $k \in \{0,1\}$, and any appropriate
index $i$. Let us prove this on the $0$-level for $\varphi $ and $\varphi'$ in $W^0(M)$.
For any $m \in M$, we have:
$$\eqalign{(d^0\varphi '\lodot d^0\varphi )(m) & = (id_M \ot \mu _H)(d^0\varphi '\ot \id _H)(d^0\varphi (m))\cr
& = (id_M \ot \mu _H)(d^0\varphi '\ot \id _H)(\varphi(m_0)_0\ot \varphi(m_0)_1 \sigma_H(m_1))\cr
& = \varphi'(\varphi(m_0)_0)_0 \ot 
\varphi'(\varphi(m_0)_0)_1\sigma_H(\varphi(m_0)_1)\varphi(m_0)_2\sigma_H(m_1)\cr
& = \varphi'(\varphi(m_0)_0)_0 \ot 
\varphi'(\varphi(m_0)_0)_1\varepsilon_H (\varphi(m_0)_1)\sigma_H(m_1) \cr
 & = (id_M \ot \mu _H)((\Delta _M \circ \varphi') \ot \id _H)[\varphi(m_0)_0 \ot 
\varepsilon_H (\varphi(m_0)_1)\sigma_H(m_1)] \cr
& = (id_M \ot \mu _H)((\Delta _M \circ \varphi')  \ot \id _H)[\varphi(m_0) \ot 
\sigma_H(m_1)] \cr
& = (id_M \ot \mu _H)((\Delta _M \circ \varphi'\circ\varphi)  \ot \sigma_H )\Delta_M(m) \cr
& = d^0(\varphi '\lodot\varphi )(m)\cr
\cr
{\hbox {and}} \qquad d^1\varphi \lodot d^1\varphi ' (m) &=  (id_M \ot \mu _H)(d^1\varphi '\ot \id _H)(d^1\varphi (m))\cr
&=  (id_M \ot \mu _H)(d^1\varphi '\ot \id _H)(\varphi (m)\ot 1)\cr
&=  (id_M \ot \mu _H)(\varphi '(\varphi (m))\ot 1\ot 1)\cr 
&=  \varphi '(\varphi (m))\ot 1\cr
& = d^1(\varphi '\lodot\varphi )(m).\cr}$$

We do not write down the computations on the $1$-level, which are very similar to the previous ones. 
We leave  to the reader the straightforward proof of $d^i(\upsilon_k) = \upsilon_{k+1}$ and also
the easy checking of the following three formulae
$$d^2d^0 = d^0d^1 , \ \ \ \ d^1d^0 = d^0d^0, \ \ \ \  d^2d^1 = d^1d^1,$$
which mean that the diagram (3) is precosimplicial.
\dm

\bigskip

We define the $0$-cohomology group $\Hr^0(H, M)$ and the $1$-cohomology set $\Hr^1(H, M)$ in the following way. Let 
$$\Hr^0(H, M) = \{ \varphi \in \Aut _S(M) \ \vert \ d^1\varphi = d^0\varphi \}$$ be the equalizer
of the pair $(d^0, d^1)$. It is obviously a group since $d^i$ is a morphism of monoids.
\goodbreak

The set $\Zr^1(H, M)$ of {\sl  $1$-Hopf cocycles of $H$ with coefficients in $M$} 
is the subset of $W_S^1(M)$ defined~by 
$$\Zr^1(H, M) = \left\{ \Phi \in W_k^1(M) \quad \left\vert 
\eqalign{
\quad & (\Zr\Cr_1) \ \ \ \Phi(ms) = \Phi(m)s {\hbox {, for all}} \ m \in M \ {\hbox {and}} \ s \in S\hfill \cr
 &(\Zr\Cr_2) \ \ \ (\id _M \ot \varepsilon _H) \circ \Phi = \id _M\hfill\cr
& (\Zr\Cr_3) \ \ \ d^2\Phi \lodot d^0\Phi = d^1\Phi\hfill\cr}\right. \right\}.$$

\noindent The group $ \Aut _S(M)$
acts on the right on $\Zr^1(H, M)$  by $$(\Phi \leftharpoonup f) = d^1f^{-1}\lodot\Phi\lodot d^0f,$$
where $\Phi \in \Zr^1(H, M)$ and $f \in  \Aut _S(M)$.
Two $1$-Hopf cocycles $\Phi$ and $\Phi'$
are said to be {\sl  cohomologous} if they belong to the same orbit
under the action of $ \Aut _S(M)$ on $\Zr^1(H, M)$. We denote by $\Hr^1(H, M)$
the quotient set $\Aut _S(M) \backslash \Zr^1(H, M)$; it 
is pointed with distinguished point the class of the map $\upsilon_1 = \id _M \ot \eta _H$.
\medskip
For $i = 0,1$, we call $\Hr^i(H, M)$ the {\sl  $i^{\rom th}$-Hopf cohomology set of $H$ with coefficients 
in the {\hbox{$(H,S)$-Hopf}} module $M$}.

%%When $H= k^G$ for a finite group 
%%$G$, we recover Serre's Galois cohomology sets $\Hr^i(G, \Aut _S(M))$, as shown in Corollary
%%3.2.

\bigskip
\noindent {\sl  1.2. The main theorem: Comparison of the $1$-Hopf cohomology set with twisted forms
in the Hopf-Galois context.}

\smallskip

\noindent Let $H$ be a Hopf-algebra, $\psi: R \lr S$ be an $H$-Hopf-Galois extension, 
and  $M = N_0 \ot _RS$  be the extended $S$-module of  an $R$-module $N_0$.
We  endow $M$ with the canonical $(H,S)$-Hopf module structure given by the coaction
$\Delta _M  = \id _{N_0} \ot \Delta _S $. 
The central result of this paper asserts that the Hopf 1-co\-ho\-mo\-logy set $\Hr^1(H, M)$ is isomorphic to the 
pointed set  ${\rm Twist }(S/R, N_0)$ of twisted forms of $N_0$ up to isomorphisms.

\medskip

Let $\psi: R \lr S$ be any extension of rings and $N_0$ be an $R$-module. Recall that a {\sl  twisted form of $N_0$
(over $S/R$)} is a pair 
$(N, \varphi )$, where $N$ is an $R$-module and $\varphi: N \ot _R S \lr N_0 \ot _R S$  is an {\hbox{ $S$-li\-near isomorphism}}. Let ${\rm twist }(S/R, N_0)$ be the set of
twisted forms of $N_0$.
Two twisted forms $(N, \varphi )$ and $(N', \varphi ')$ of $N_0$ are
{\sl  isomorphic} if $N$ and $N'$ are isomorphic as $R$-modules.
Following [6], we denote by ${\rm Twist }(S/R, N_0)$ the pointed set of
isomorphism classes of twisted forms of $N_0$, the distinguished point being the 
class of $(N_0, {\rm id}_{N_0}\ot {\rm id}_{S})$. 
We mention here that all the  results of [10] involving equivalence classes of twisted forms are actually proven for this definition
of ${\rm Twist }(S/R, N_0)$  and not for the one given in [10, \pa 6.3], where the  equivalence relation is too restrictive.

\medskip
\noindent {\bf Theorem 1.2.} {\sl  Let $H$ be a Hopf-algebra, $\psi: R \lr S$ be an $H$-Hopf-Galois extension, 
and $M = N_0 \ot _RS$ be the extended $S$-module of an $R$-module $N_0$. There is an isomorphism of pointed sets
$$ \Hr^1(H, M) \cong {\rm Twist }(S/R, N_0).$$ }

Theorem 1.2 allows us to state  the following
noncommutative generalization of 
Noether's cohomological form of Hilbert's Theorem 90. 

\medskip
\noindent {\bf Corollary 1.3.} {\sl  Let $H$ be a Hopf-algebra and $\psi: K \lr L$ be an $H$-Hopf-Galois extension of
 division algebras. Then, for any positive integer $n$, we have $$\Hr^1(H, L^n) = \{1\}.$$
}

\noindent Here we denote by $1$ the distinguished point of $\Hr^1(H, L^n)$. 
\medskip

\noindent {\sl  Proof of Corollary 1.3.} Observe that 
$L^n$ is isomorphic to the extended $L$-module $K^n \ot _KL$. By Theorem 1.2, the pointed set $\Hr^1(H, L^n)$ is isomorphic to ${\rm Twist }(L/K, K^n)$, which is known to be trivial
([10, Corollary 6.21]). \dm

\medskip

The rest of the paper is mainly devoted to the proof of Theorem 1.2. This is done
in two steps. At first we introduce a non-abelian cohomology theory $\Dr^i(H, M)$, for $i=0,1$,
which is  related to noncommutative descent theory.
In Theorem 2.6, we prove  the isomorphism $ \Dr^1(H, M) \cong {\rm Twist }(S/R, N_0)$. 
Subsequently we show that the Hopf cohomology sets $\Hr^i(H, M)$ are isomorphic to the
descent cohomology sets $\Dr^i(H, M)$.

\bigskip
\goodbreak
%\vfill
 
\noindent {\bf 2.  Descent cohomology sets.}

\smallskip
\noindent In this section we introduce two descent cohomology sets. We compute 
them in the Galois case and relate them to the usual non-abelian group cohomology theory.
In addition, in the Hopf-Galois context, we prove that the $1$-descent cohomology set
classifies twisted forms  and interpret it
in terms of torsors on the module of coefficients.

\smallskip

\noindent {\sl  2.1. Definition of descent cohomology sets.}

\smallskip

\noindent Let $H$ be a Hopf-algebra, $S$ be an $H$-comodule algebra,
and $M$ be an $(H,S)$-Hopf module  
with coaction $\Delta _M:  M \lr M \ot H$. 
We define the $0$-cohomology group $\Dr^0(H, M)$ %%and a $1$-cohomology set $\Dr^1(H, M)$ in the following way. Let 
by $$\Dr^0(H, M) = \{ \alpha \in \Aut _S(M) \ \vert \ (\alpha\ot \id_H) \circ \Delta _M = \Delta _M \circ \alpha \}.$$
It is the set of the $S$-linear automorphisms of $M$ which are maps of $H$-comodules. This set obviously carries
a group structure given by the composition of automorphisms.

\medskip
\noindent {\bf Lemma 2.1.}
{\sl   Let $H$ be a Hopf-algebra and $S$ be an $H$-comodule algebra.
Any isomorphism {\hbox{$f : M \to M'$}}
of $(H,S)$-Hopf modules induces an isomorphism of groups $f^*: \Dr^0(H, M') \lr \Dr^0(H, M)$
given on $\alpha' \in \Dr^0(H, M')$ by:
$$f^*\alpha' = f^{-1}\circ \alpha' \circ f.$$

}

}

\medskip
\Dem 
The  $S$-linearity of $f^*\alpha'$ immediately follows from
the $S$-linearity of $f$ and that of $\alpha'$.
In order to prove that $f^*\alpha'$ belongs to $\Dr^0(H, M)$, it is sufficient
to observe that the following 
diagram is commutative.
\medskip

%%DIAG
  $$\xymatrix{  	M\ar [ddd]_{\Delta _M} \ar [rrrr]^{f^*\alpha'}\ar[dr]_{f}&&&&M
\ar [ddd]^{\Delta _M}\ar[dl]^{f}  \\
   			&M'\ar[rr]^{\alpha'}\ar[d]_{\Delta _{M'}}&&M'\ar[d]^{\Delta _{M'}}& \\
			&M'\ot H\ar[rr]_{\alpha' \otimes {\rm id}_H} &&M'\ot H&\\
			M\ot H \ar [rrrr]_{f^*\alpha'\otimes {\rm id}_H}\ar[ru]^{f\otimes{\rm id}_H}&&&&
M\ot H\ar[ul]_{f\otimes {\rm id}_H}\\}$$
\dm

\goodbreak

We introduce now a $1$-cohomology set $\Dr^1(H, M)$ in the following way.
The set $\Cr^1(H, M)$ of {\sl  $1$-descent cocycles of $H$ with coefficients in $M$} is defined to be
 the set of all $k$-linear $H$-coactions $F: M \lr M \ot H$ on $M$
making $M$ an $(H,S)$-Hopf module. In other words, one has:
$$\Cr^1(H, M) = \left\{ F: M \lr M \ot H \quad \left\vert 
\eqalign{
\quad & (\Cr\Cr_1) \ \ \ F(ms) = F(m)\Delta_S(s)  {\hbox {, for all}} \ m \in M \ {\hbox {and}} \ s \in S\cr
 &(\Cr\Cr_2) \ \ \ (\id _M \ot \eps _H) \circ F = \id _M\cr
& (\Cr\Cr_3) \ \ \ (F \ot \id _H) \circ F = (\id _M \ot \Delta _H) \circ F\cr}\right. \right\}.$$

\noindent Notice that $\Cr^1(H, M)$ is pointed (hence not empty) with the  coaction
map $\Delta _M$ as distinguished point.
%% and that  $\Cr^1(H, M)$ is exactly the set of all descent data on $M$, as these are defined in [11]. 
\goodbreak
\medskip

\noindent {\bf Lemma 2.2.}
{\sl   Let $H$ be a Hopf-algebra and $S$ be an $H$-comodule algebra.
Any isomorphism
{\hbox{$f : M \to M'$}} of $S$-modules induces a bijection
 $f^*: \Cr^1(H, M') \lr \Cr^1(H, M)$ given on $F' \in \Cr^1(H, M')$ by 
 $$f^*F' = (f^{-1}\ot \id_H) \circ F \circ f.$$
For any $S$-module $M$, one has $$(\id _{M})^* = \id _{{\rm C}^1(H, M)}.$$ For any
composable isomorphisms of $S$-modules $f: M \lr M'$ and $f': M' \lr M''$, the following
equality holds $$(f'\circ f)^* = f^* \circ f'^*.$$

\noindent If moreover $f: M \lr M'$ is an isomorphism
of $(H,S)$-Hopf modules, then $f^*$ realizes an isomorphism of
pointed sets between $\Cr^1(H, M')$ and $\Cr^1(H, M)$.

}

\medskip
\Dem Let  $f: M \lr M'$ be an isomorphism  of $S$-modules. 
The  $(H,S)$-linearity of $f^*F'$ immediately follows from
the $S$-linearity of $f$ and from the $(H,S)$-linearity of $F'$. The coassociativity of
$f^*F'$
comes from the commutativity of the
diagram
%%DIAG
  $$\xymatrix{  	M\ar [ddd]_{f^*F'} \ar [rrrr]^{f^*F'}\ar[dr]_f&&&&M\ot H
\ar [ddd]^{f^*F'\otimes{\rm id}_H}\ar[dl]^{f\otimes{\rm id}}  \\
   			&M'\ar[rr]^{F'}\ar[d]_{F'}&&M'\ot H\ar[d]^{F'\otimes 1}& \\
			&M'\ot H\ar[rr]_{{\rm id}_{M'}\otimes\Delta_H} &&M'\ot H\ot H&\\
			M\ot H \ar [rrrr]_{{\rm id}_M\otimes\Delta_H}\ar[ru]^{f\otimes{\rm id}_H}&&&&
M\ot H\ot H\ar[ul]_{\ f\otimes {\rm id}_H^{\otimes 2}}\\} $$
\noindent whereas the compatibility of $f^*F'$ with the counity of $H$
is expressed by the commutativity of the diagram
$$
\xymatrix{  	M  \ar [rrr]^{f^*F'}\ar[dr]_f&&&M\ot H \ar[dl]^{f\otimes{\rm id}} \\
   			&M'\ar[dr]_{{\rm id}_{M'}}\ar[r]^{F'} &M'\ot H\ar[d]^{{\rm id}_{M'}\otimes \varepsilon_H}& \\
			&  &M'\ &\\
			 &&&M\ar [uuu]_{{\rm id}_{M}\otimes\varepsilon_H}\ar[ul]^{f }\\}$$
Hence we have shown that  $f^*F'$ belongs to $\Cr^1(H, M)$.
By the very definition, $f^*F'$ is bijective and $(\id _{M})^* = \id _{{\rm C}^1(H, M)}$.

\medskip
Let $f: M \lr M'$ and $f': M' \lr M''$ be two isomorphisms of $S$-modules.
One has, for any $F' \in \Cr^1(H, M')$, the following equalities
$$(f'\circ f)^* (F') = \bigl((f'\circ f)^{-1}\ot \id_H\bigr) \circ F'  \circ (f'\circ f)
= \bigl((f^{-1}\circ f'^{-1})\ot \id_H\bigr) \circ F'  \circ (f'\circ f) = f^*(f'^*F').$$

Moreover, if $f $ is an isomorphism
of $(H,S)$-Hopf modules, the map $f^*$ preserves the distinguished points: indeed,
the equality $f^*\Delta _{M'} = \Delta _M$
is equivalent to the fact that $f$ is a morphism of $(H,S)$-Hopf modules.
\dm

\medskip
From Lemma 2.2, one readily obtains the following result:
\medskip

\noindent {\bf Corollary 2.3.}
{\sl  Let $H$ be a Hopf-algebra, $S$ be an $H$-comodule algebra,
and $M$ be an $(H,S)$-Hopf module. The group $ \Aut _S(M)$
acts on the right on $\Cr^1(H, M)$  by $$(F \leftharpoonup f) = f^*F = (f^{-1}\ot \id_H) \circ F  \circ f,$$
where $F \in \Cr^1(H, M)$ and $f \in  \Aut _S(M)$.}

\medskip

\noindent Two $1$-descent cocycles $F$ and $F'$
are said to be {\sl  cohomologous} if they belong to the same orbit under
 the action of $ \Aut _S(M)$ on $\Cr^1(H, M)$. We denote by $\Dr^1(H, M)$
the quotient set $\Aut _S(M) \backslash \Cr^1(H, M)$; it 
is pointed with distinguished point the class of the coaction $\Delta _M$.

\medskip

For $i = 0,1$, we call $\Dr^i(H, M)$ the {\sl  $i^{\rom th}$-descent cohomology set of $H$ with coefficients 
in $M$}. The choice of this name finds its motivation in the following observation. Suppose that $\psi: R \lr S$ is an
$H$-Hopf-Galois extension. As shown in [11], an $(H,S)$-Hopf module may always be descended to an $R$-module $N_0$, that is
$M$ is isomorphic to an extended  $S$-module $ N_0 \ot _RS$. 
The set  $\Cr^1(H, M)$ is exactly those  of all descent data on $M$ described in  [10].

\medskip
\noindent {\bf Corollary 2.4.}
{\sl  Let $H$ be a Hopf-algebra and $S$ be an $H$-comodule algebra. 

\item{--} Any isomorphism
$f: M \lr M'$ of $S$-modules induces a bijection $f^*: \Dr^1(H, M') \lr \Dr^1(H, M)$. 

\item{--} Any isomorphism
{\hbox{$f: M \lr M'$}} of $(H,S)$-Hopf modules induces an isomorphism of pointed sets $f^*:
\Dr^1(H, M') \lr \Dr^1(H, M)$.

}

\medskip
\Dem Suppose that $F_1$ and  $F_2$ are two cohomologous 1-cocycles of $\Cr^1(H, M')$, with $g \in \Aut _S(M')$  such that
$ F_1 = g^*F_2$. Then $f^*F_2 = f^*g^*F_1 = f^*g^*(f^{-1})^*f^*F_1 = (f^{-1}gf)^*(f^*F_1)$, so
 $f^*F_1$ and $f^*F_2$ are cohomologous in $\Cr^1(H, M)$.
\dm

\bigskip
\noindent {\sl  2.2. Application to the Galois case.}

\smallskip

\noindent We work now with the Hopf algebra $k^G$  dual  to the group algebra $k[G]$ for $G$ a finite group. Let 
 $\psi: R \lr S$ be a $k^G$-Galois extension  and $M$    a $(G,S)$-Galois module.
We may assume that $M$ is already extended, so that $M$ is equal to $N_0 \ot _RS$ for  an $R$-module  $N_0$.
Endow $M$ with the canonical $(H,S)$-Hopf module structure given by the coaction
$\Delta _M  = \id _{N_{0}} \ot \Delta _S $.
In this paragraph, we compute the descent cohomology set of $k^G$ with coefficients 
in  $M = N_0 \ot _RS$ in terms of the Galois $1$-cohomology set of $G$ with coefficients in $\Aut _S(M)$.
\medskip

Recall that for any group $G$ and any
(left) $G$-group $A$, one classically defines  two non-abelian cohomology sets 
of $G$ with coefficients in $A$ (see [12] and [13]). This is done in the
following way.
The 0-cohomology group
$\Hr^0(G,A) $ is the group $A^G$ of 
invariant elements of $A$ under the action of $G$.
The set $\Zr^1(G, A)$ of 1-cocycles is given by
$$\Zr^1(G, A) = \{ \alpha \in {\goth Set}(G, A) \ \vert \ \ \alpha (gg') =
\alpha (g){{}^{g}\! \bigl(}\alpha (g')\bigr),
 \ \ \forall \ g, g' \in G \}.$$
It is pointed with distinguished point the constant map $1: G \lr A$.

\goodbreak
The group
$A$ acts on the right on $\Zr^1(G, A)$  by $$(\alpha \leftharpoonup a)(g) = a^{-1}\alpha (g) \ {}^{g}\! a,$$
where $a \in A$, $\alpha  \in \Zr^1(G, A)$, and $g \in G$.
Two 1-cocycles $\alpha $ and $\alpha '$
are {\sl  cohomologous} if they belong to the same orbit under this action.
The non-abelian 1-cohomology set  $\Hr^1(G, A)$ is the left quotient
$A \backslash \Zr^1(G, A)$.
Then $\Hr^1(G,A)$ is pointed with distinguished point the class of the constant  {\hbox{map  $1: G \lr A$.}}

\medskip Let $G$ be a finite group, $\psi: R \lr S$ be a $G$-Galois extension, 
and $M = N_0 \ot _RS$  be the extended $S$-module of  an $R$-module $N_0$.
The $S$-module $M$ is a $(G,S)$-Galois  module by the canonical action given on an indecomposable tensor $n \ot s \in N_0 \ot _RS$ by
$$g(n \ot s) =n \ot g(s),$$
where $g \in G$, $n \in N_0$, and $s \in S$.
The group
$G$ acts by automorphisms on  ${\rm Aut}_S(M)$ by 
$${^g\!f} = ({\rm id}_{N_0}\ot g)\circ f \circ ({\rm id}_{N_0}\ot g^{-1}),$$
where $g\in G$ and $f \in {\rm Aut}_S(M)$. Hence $\Aut _S(M)$ becomes a $G$-group and we get at our disposal the two non-abelian cohomology
sets $\Hr^0(G, \Aut _S(M))$ and $\Hr^1(G, \Aut _S(M))$.

\medskip
\noindent {\bf Proposition 2.5.} {\sl  Let $G$ be a finite group, $\psi: R \lr S$ be a $G$-Galois extension, 
and $M = N_0 \ot _RS$  be the extended $S$-module of  an $R$-module $N_0$. There is the   equality of groups
$$\Dr^0(k^G, M) = \Hr^0(G, \Aut _S(M))$$ and an isomorphism of pointed sets
$$ \Dr^1(k^G, M) \cong \Hr^1(G, \Aut _S(M)).$$ }

\Dem
Let us prove the equality between the groups. It is sufficient to show that  for any  $f \in \Aut _S(M)$,
the condition  $(f\ot \id_{k^G}) \circ  \Delta _M  = \Delta _M \circ f$ is equivalent to the fact that $f$ is $G$-invariant. Indeed, the first condition   
reflects that $f$ belongs to $\Dr^0(k^G, M)$, whereas   $\Hr^0(G, \Aut _S(M))$ is precisely the group
$\Aut _S(M)^G$ of $G$-invariant automorphisms in $\Aut _S(M)$. Pick
$f \in \Aut _S(M)$,  $n\in N_{0}$, and $s\in S$. One has   
$$\bigl((f\ot \id_{k^G}) \circ  \Delta _M \bigr)(n\ot s) = \sum _{g\in G}(f\ot \id_{k^G})\bigl(n\ot g(s) \ot \delta _g\bigr) =
\sum _{g\in G}\bigl(f\circ (\id_{N_0} \ot g)\bigr)(n\ot s) \ot \delta _g.$$
On the other hand, setting  $f (n \ot s) = n' \ot s'$, one gets
$$(\Delta _M \circ f) (n \ot s) = \Delta _M (n' \ot s') 
= \sum _{g\in G}  \bigl(n'\ot g(s')\bigr) \ot \delta _g 
= \sum _{g\in G}\bigl((\id_{N_0} \ot g)\circ f\bigr)(n\ot s) \ot \delta _g.
$$
Since $\{ \delta _g \}_{g \in G}$ is a basis of $k^G$, the relation $(f\ot \id_{k^G}) \circ  \Delta _M  = \Delta _M \circ f$  
is equivalent to the set of equalities $f\circ (\id_{N_0} \ot g) = (\id_{N_0} \ot g)\circ f$, with  $g$ running through $G$. This exactly means
that $f$ is $G$-invariant in $\Aut _S(M)$. 
\medskip
\goodbreak
We prove now the isomorphism on the 1-cohomology level. 
Let us show that  any $F \in \Cr^1(k^G, M)$ induces a $(G,S)$-Galois module 
action $\gamma \in \Aut^\gamma _{S}(M)$
 defined  by
$$F(m) = \sum _{g \in G}\bigl(\gamma (g)\bigr)(m) \ot \delta _g.$$
For simplicity  denote $\gamma (g)(m)$ by $g(m)$. The $k$-linearity of $F$ tells us that $g(m + m') = g(m) + g(m')$, for any $g \in G$ and 
$m, m' \in M$; the equality $(\id _M \ot \eps _{k^G}) \circ F = \id _M$ implies that $1(m) = m$;
the coassociativity condition of $F$ says  that $(gg')(m) = g\bigl(g'(m)\bigr)$, for any $g,g' \in G$ and 
$m\in M$;
finally the $(k^G,S)$-linearity of $F$ is  equivalent to the $(G,S)$-linearity of $\gamma$.
As shown in [10], the 	action map $\gamma$ gives rise to the 1-Galois cocycle $\alpha: G \lr \Aut _S(M)$ defined by
$$\alpha (g) = \gamma (g)\circ (\id_{N_0}\ot  g^{-1}).$$ It is easy to check that the correspondence between $F$
and $\alpha$ is bijective. Thus already at the 1-cocycle level there exists  a bijection between $\Zr^1(G, \Aut _S(M))$ and $ \Cr^1(k^G, M)$.
  %%: for  any 1-Galois cocycle $\alpha: G \lr \Aut _S(M)$, one defines 
%%$\gamma (g) = \alpha (g)\circ (\id_{N_0}\ot  g)$ and $F(m) = \sum _{g \in G}\bigl(\gamma (g)\bigr)(m) \ot \delta _g$

Take two cocycles $F$ and $F'$ in $\Cr^1(k^G, M)$. Denote by $\gamma $ (respectively $\gamma '$) the corresponding
Galois actions and by $\alpha$ (respectively $\alpha '$) the Galois cocycles associated with $\gamma $ (respectively $\gamma '$).
Suppose that the cocycles $F$ and  $F'$ are cohomologous, with $f \in \Aut _S(M)$  such that
$(f\ot \id_{k^G}) \circ F = F' \circ f $. Then $f \circ \gamma (g) =
\gamma '(g) \circ f$, for all $g \in G$, or equivalently $\gamma (g) = f^{-1} \circ \gamma' (g) \circ f$. Therefore
 $$\eqalign{\alpha (g) & = f^{-1} \circ \gamma '(g) \circ f \circ (\id _{N_0} \ot g^{-1}) \cr & = 
f^{-1} \circ \gamma '(g) \circ (\id _{N_0} \ot g^{-1}) \circ (\id _{N_0} \ot g)\circ f \circ (\id _{N_0} \ot g^{-1})\cr
&= f^{-1} \circ \alpha '(g) \circ {^g\! f},\cr}$$
which means that $\alpha$  and $\alpha '$ are Galois-cohomologous. 
Conversely, the previous equalities show that 
two cohomologous  Galois cocycles $\alpha$  and $\alpha '$ give rise to two cohomologous cocycles  $F$ and $F'$ in $\Cr^1(k^G, M)$.
 \dm

\bigskip
\goodbreak
\noindent {\sl  2.3. Comparison between the $1$-descent cohomology set and the set of twisted forms
in the Hopf-Galois context.}

\smallskip

\noindent Let $H$ be a Hopf-algebra, $\psi: R \lr S$ be an $H$-Hopf-Galois extension, 
and $M = N_0 \ot _RS$  be the extended $S$-module of  an $R$-module $N_0$.
We  endow $M$ with the canonical $(H,S)$-Hopf module structure given by the coaction
$\Delta _M  = \id _{N_0} \ot \Delta _S $. 
The main result of this paragraph asserts that the descent {\hbox {1-coho\-mology}} set $\Dr^1(H, M)$ is isomorphic to the 
pointed set  ${\rm Twist }(S/R, N_0)$ of twisted forms of $N_0$ up to isomorphisms.

\goodbreak

\medskip
\noindent {\bf Theorem 2.6.} {\sl  Let $H$ be a Hopf-algebra, $\psi: R \lr S$ be an $H$-Hopf-Galois extension, 
and $M = N_0 \ot _RS$ be the extended $S$-module of an $R$-module $N_0$. Then 
there is an isomorphism of pointed sets
$$ \Dr^1(H, M) \cong {\rm Twist }(S/R, N_0).$$ }

In order to prove Theorem 2.6, we need an intermediate result.  For any   $F \in \Cr^1(H, M)$ denote by $N_F$  the $R$-module of $F$-coinvariants, that is
$N_F = \{ m \in M \ \ \vert \ \ F(m) = m\ot 1\}$. We state the following lemma:

\medskip

\goodbreak
\noindent {\bf Lemma 2.7.} {\sl  Under the same hypotheses as in Theorem 2.6, 
for any   $F \in \Cr^1(H, M)$, there exists an isomorphism 
$$\varphi _F: N_F \ot _R S {\buildrel \sim \over \lr}  M$$ given by $\varphi _F (m \ot s) = ms$, for any
$m \in N_F$ and $s \in S$. }

\medskip
\noindent {\sl  Proof.} 
The existence of the isomorphism $\varphi _F $ results from Hopf-Galois descent theory [11, Theorem 3.7] (see 	also [5]). Indeed, consider the
functor ``{\sl  restriction of scalars}" $\psi ^*: {\goth Mod}_S \lr {\goth Mod}_R$  and its  
left adjoint functor $\psi_!: {\goth Mod}_R \lr {\goth Mod}_S$, the functor  ``{\sl  extension of scalars}". Then
$\varphi _F $ is nothing but a counit 
for the comonad on ${\goth Mod}_S$ induced by the adjunction 
$\psi _{!}\ \adj \ \psi ^*$  (see, {\sl  e.g.}, [4]). 

We explicit now the expression of $\varphi _F $. By arguments stemming from descent theory ([3], [10]), the $S$-module $M$ is isomorphic to $N_d \ot _RS$, where 
$N_d$ is the $R$-module deduced from the Cipolla descent data $d$ on $M$ associated to
the $(H,S)$-Hopf module structure of $M$. By [10, Prop. 4.10],
$d$ is the map given by the composition
$$\matrix {
M&\hfll{F}{}& M \ot H &\hfll{\beta}{}&M\ot _S(S\ot H ) &
\hfll{{\ninerom id} _M \otimes \Gamma _\psi^{-1}}{}&M\ot _S(S\ot _RS)&\hfll{\beta '}{}&
M\ot _RS,\cr}$$ 
where  $\beta $ (respectively $\beta '$) is the obvious $k$-linear 
(respectively $S$-linear) isomorphism and $\Gamma _\psi$ is the Galois isomorphism mentioned in the Conventions.

Let us now compute $d$. For $m \in M$, set $F(m) = \sum_im_i \ot h_i \in M \ot H$.
For any fixed index $i$, set $\Gamma _\psi^{-1}(1 \ot h_i) = \sum _j s_{ij} \ot t_{ij}$, or equivalently
$\sum _j s_{ij} \Delta_S(t_{ij})$ $ =  1 \ot h_i$.
So $$d(m) = \sum_i\sum_jm_is_{ij} \ot t_{ij}.$$
According to [10, Cor. 4.11], we have
$N_d = \{ m \in M \ \ \vert \ \ \sum_i\sum_jm_is_{ij} \Delta_S (t_{ij}) = m \ot 1\}$, therefore
$$N_d = 
 \{ m \in M \ \ \vert \ \ \sum_im_i\ot h_i = m \ot 1\} =  \{ m \in M \ \ \vert \ \ F(m) = m \ot 1\}  = N_F.$$
It is proven in [3] that the descent isomorphism from $N_d \ot _RS$ to $M$ 
is given by the correspondence $m \ot s \longmapsto ms$
for 
$m \in N_d$ and $s \in S$. %%%%%%
\dm

\medskip
\noindent {\sl  Proof of Theorem 2.6.} Let $F$ be an element of $\Cr^1(H, M)$ and $\varphi _F$ be the isomorphism
from $N_F \ot _RS$ to $M = N_0\ot _RS$ given by the previous lemma. The datum  $(N_F, \varphi _F)$ is a twisted form  of
$N_0$. Denote by $\tilde {\T}$ the map from $\Cr^1(H, M)$ to the set ${\rm twist }(S/R, N_0)$ 
 defined by
$$\tilde {\T} (F) = (N_F, \varphi _F).$$ The map $\tilde {\T}$ obviously sends the distinguished point $\Delta _M$
of $\Cr^1(H, M)$ to the distinguished point $(N_0, \id _{N_0\ot _RS})$ of ${\rm twist }(S/R, N_0)$.

Suppose that $F$ and $F'$ are cohomologous in $ \Cr^1(H, M)$. We claim that the corresponding descended modules
$N_F$ and $N_{F'}$ are isomorphic in ${\goth Mod}_R$. Indeed, let 
  $f \in  \Aut _S(M)$ such that $(f\ot \id_{H}) \circ F = F' \circ f$. 
For any $n \in N_F$, the image $f(n)$ belongs to $N_{F'}$,
since $$F'(f(n)) = (f\ot \id_H) (F(n)) = (f\ot \id_H) (n \ot 1) = f(n) \ot 1.$$ 
So the automorphism $f$ induces an isomorphism 
from $N_F$ to $N_{F'}$. From this fact we deduce a quotient map 
$$\T: \Dr^1(H, M) \lr {\rm Twist }(S/R, N_0).$$

We now prove that $\T $ is an isomorphism of pointed sets. In order to do this, we introduce the map {\hbox{$\tilde {\D}: {\rm twist }(S/R, N_0) \lr \Cr^1(H, M)$}}
which associates to any twisted form $(N, \varphi)$ of $M$ 
the
map $F_N: M \lr M \ot H$ defined
by $$F_N = (\varphi ^{-1})^* (\id _N \ot \Delta _S) = (\varphi \ot \id _H) \circ (\id _N \ot \Delta _S) \circ \varphi^{-1}.$$
Since $(\id _N \ot \Delta _S)$ is the canonical  $(H,S)$-Hopf module structure on $N \ot _RS$, by Lemma 2.2,
the map $F_N$ belongs to $\Cr^1(H, M)$.

Suppose that $(N, \varphi)$ and $(N', \varphi')$ are two equivalent twisted forms of $M$ via $\vartheta \in \Aut _S(M)$. 
Set $f = \varphi ' \circ (\vartheta \ot \id _S) \circ \varphi ^{-1}$. Observe that  the following diagram commutes:

   $$\xymatrix{ M \ar[rr]^{\varphi ^{-1} }\ar[dd]_{f}
\ar@/ ^2.1pc/[rrrrrr]^{F_N}&& N \ot _RS \ar[rr]^{{\rm id}_N \otimes \Delta _S  }
\ar[dd]^{\vartheta \otimes {\rm id}_S} && N \ot _RS \ot H
\ar[rr]^{\varphi \otimes {\rm id}_H}
\ar[dd]_{\vartheta \otimes {\rm id}_S\otimes {\rm id}_H} && M \ot H\ar[dd]^{f\otimes {\rm id}_H} \\
\\
M \ar[rr]_{\varphi '^{-1} }\ar@/ _2.1pc/[rrrrrr]_{F_{N'}}&& N' \ot _RS \ar[rr]_{{\rm id}_{N'} \otimes \Delta _S  } && N' \ot _RS \ot H
\ar[rr]_{\varphi' \otimes {\rm id}_H} && M \ot H      }$$

\smallskip

\noindent So $F_{N'}$ equals $f^*F_N$ and therefore $\tilde {\D}$ induces a quotient map
$$\D:   {\rm Twist }(S/R, N_0)\lr  \Dr^1(H, M).$$ 

It remains to prove that $\T \circ \D$ and $\D \circ \T$ are the identity maps.

\medskip
\noindent {\sl  The composition $\T \circ \D$ is the identity.} Let $(N, \varphi)$ be a twisted form of $N_0$. Since $N_{F_N}\ot _RS $ is isomorphic to
$N\ot _RS $ (Lemma 2.7), we deduce from Hopf-Galois descent theory [11, Theorem 3.7] the existence of an isomorphism 
$\vartheta: N \lr N_{F_N}$. So the twisted form $\tilde {\T} (\tilde {\D} (N, \varphi)) $ is equivalent to $(N, \varphi)$.
In concrete terms, $\vartheta$ fits into the   following commutative diagram of $R$-modules with exact rows:
$$\xymatrix{ 	0\ar[r]	& N \ar[r] \ar[d]^\vartheta _\wr& N\ot_RS\ar[d]^\varphi _\wr \ar@<0.5ex>[rr]^{{\rm id}_N\otimes \Delta_S \ } 
\ar@<-0.5ex>[rr]_{({\rm id}_{N\otimes_R S})\otimes\eta_H \ }  && N\ot_RS\ot H\ar[d]^{\varphi \otimes {\rm id}_H}_\wr 	\\ 
     0\ar[r]& N_{F_N\ } \ar@{^{(}->}[r]  & M  \ar@<0.5ex>[rr]^{F_N \ } \ar@<-0.5ex>[rr]_{{\rm id}_{M}\otimes\eta_H \ } 
&& M \ot H\\ }
$$
Hence one gets \ $\T \circ \D = \id$.

\goodbreak
\medskip
\medskip
\noindent {\sl  The composition $\D \circ \T$ is the identity.}
Let $F$ be an element of $\Cr^1(M, H)$. Consider the following diagram:
\smallskip
  $$\xymatrix{ M\ar[rrrr]^{F_{N_F}}\ar[rd]^{\varphi _F^{-1}}\ar@{=}[dd]  &   & &    &  
M \ot H \ar[ld]_{\varphi _F^{-1}\otimes{\rm id}_H}\ar@{=}[dd]\\
    & N_F \ot _RS \ar[rr]^{{\rm id}_{N_F} \otimes \Delta _S}\ar[ld]^{\varphi _F}&& 
N_F \ot _RS\ot H\ar[rd]_{\varphi _F\otimes{\rm id}_H} && \\
M  \ar[rrrr]_{F}  &   & &    &  M \ot H }$$
\smallskip

\noindent The left and right triangles are trivially commutative. The upper trapezium commutes
by the definition of $F_{N_F}$. Let us show  the commutativity of the lower trapezium. 
Pick an indecomposable tensor $m \ot s$  in $N_{F}\ot _RS$. Setting $\Delta _S(s) = s_0 \ot s_1$, we have
$$(\varphi _F\ot \id _H) \circ (\id _{N_F} \ot \Delta _S)(m \ot s) =  \varphi _F(m \ot s_0) \ot s_1
=  ms_0 \ot s_1.$$ The latter equality comes from Lemma 2.7.
On the other hand, using  the $(H,S)$-linearity of $F$, one has $$(F\circ \varphi _F) (m\ot s) = F(ms) = F(m)\Delta_S(s) =  ms_0 \ot s_1.$$ 
So the whole diagram is commutative. Hence we obtain $F = F_{N_F}$,
which means $\tilde{\D} \circ \tilde{\T} = \id$. Therefore we conclude $\D \circ \T = \id$. \dm
\medskip

\bigskip
\goodbreak
\noindent {\sl  2.4. The $1$-descent cohomology set and torsors.}

\smallskip

\noindent Let $G$ be a finite group and $A$ be a $G$-group. Recall that 
an  {\sl  $A$-torsor} (or 
{\sl   $A$-principal homogeneous space}) is a non-empty $G$-set $P$ on which $A$ acts on the
right in a compatible way with the $G$-action and such that $P$ is an affine space over $A$ (see [13]).
Pursuing our analogy between non-abelian group- and Hopf-cohomology theories, we are led to state the following definition.

\goodbreak
Let $H$ be a Hopf algebra and $M$ be an $(H,S)$-Hopf module. An {\sl  $M$-torsor} is a triple
$(X, \Delta _X, \beta)$, where $\Delta _X: X \lr X \ot H$ is a map conferring $X$ 
a structure of {\hbox{$(H,S)$-Hopf}} module and $\beta: M \lr X$ is an $S$-linear isomorphism.
Denote by $\tors (M)$ the set of $M$-torsors. It is pointed with distinguished point
$(M, \Delta _M, \id _M)$.
We say that two $M$-torsors $(X, \Delta _X, \beta)$ and $(X', \Delta _{X'}, \beta')$ are {\sl  equivalent}
if there  exists  $f \in \Aut _S(M)$ such that the composition $\beta \circ f \circ \beta '^{-1}: X' \lr X$ is a morphism
of {\hbox{$(H,S)$-Hopf}} modules.
Denote by $\Tors (M)$ the set of equivalence classes of $M$-torsors; it is pointed with distinguished point
the class of $(M, \Delta _M, \id _M)$. We have the following result:

\medskip
\noindent {\bf Proposition  2.8.} {\sl  Let $H$ be a Hopf algebra and $M$ be an $(H,S)$-Hopf module.
There is an isomor\-phism of pointed sets
$$ \Dr^1(H, M) \cong \Tors(M).$$}

\Dem Define $\tilde {\U}: \Cr^1(H,M) \lr \tors (M)$ and $\tilde {\V}: \tors (M) \lr \Cr^1(H,M)$ by 
$$\tilde {\U}: F \longmapsto (M, F, \id _M) \quad {\hbox {\rom and}} \quad \tilde {\V}: 
(X, \Delta _X, \beta) \longmapsto \beta ^*\Delta _X.$$ We set here
$\beta ^*\Delta _X = (\beta^{-1}\ot \id_H) \circ \Delta _X \circ \beta$, which,  following Lemma 2.2, is an element of $\Cr^1(H,M)$ since 
$\Delta _X$ belongs to $\Cr^1(H,X)$. Using again Lemma 2.2, it is easy to check that $\tilde {\U} $ and $\tilde {\V}$ define maps 
$\U: \Dr^1(H,M) \lr \Tors (M)$ and $\V: \Tors (M) \lr \Dr^1(H,M)$ on the quotients.

It is straightforward to prove $\tilde {\V}\circ \tilde {\U} = \id _{\Cri^1(H,M)}$. Moreover, 
the torsor $(\tilde {\U}\circ \tilde {\V})(X, \Delta _X, \beta)$ equals $(M, \beta ^*\Delta _X, \id_M)$,
which, via $f = \id _M$, is equivalent to $(X, \Delta _X, \beta)$ in $\tors (M)$. \dm
\bigskip

\noindent {\bf 3. The isomorphism between Hopf cohomology sets and descent cohomology sets.}

\smallskip
\noindent In this paragraph, we interpret the noncommutative Hopf cohomology sets  in terms of the descent cohomology sets.

\smallskip
Let $H$ be a Hopf-algebra
and $(M, \Delta _M:  M \lr M \ot H)$ be an $H$-comodule.
We define a map $\tilde \kappa$ from $\Wr^1_k(M)$ to itself  by the formula
$$\tilde \kappa (\Phi) = \Phi \lodot \Delta _M,$$ 
for any $\Phi \in \Wr^1_k(M)$. The map $\tilde \kappa$ is a bijection.
Indeed, denote by $\Delta _M'$ the map $(\id_M \ot \sigma _H)\circ \Delta _M$, which is easily seen to be
the $\lodot$-inverse of $\Delta _M$ in $\Wr^1_k(M)$.
The inverse map of $\tilde \kappa$ is therefore given by
$$\tilde \kappa^{-1} (\Phi) = \Phi \lodot \Delta _M'.$$

\medskip
\noindent {\bf Theorem 3.1.} {\sl  Let $H$ be a Hopf-algebra, $S$ be an $H$-comodule algebra,
and $M$ be an $(H,S)$-Hopf module  
with coaction $\Delta _M:  M \lr M \ot H$. The identity map $\id _{{\rm Aut}_S(M)}$ realizes the equality
of groups
$$\Hr^0(H, M) = \Dr^0(H, M).$$ The translation
map $\tilde \kappa$  induces an isomorphism of pointed sets 
$$\kappa: \Hr^1(H, M) \lr \Dr^1(H, M).$$

}

\noindent As a consequence of this result and of Proposition 2.5, one immediately gets the following corollary
which relates non-abelian Hopf-cohomology objects to non-abelian group-cohomology objects:

\medskip
\noindent {\bf Corollary 3.2.} {\sl  Let $G$ be a finite group, $\psi: R \lr S$ be a $G$-Galois extension, 
and $M = N_0 \ot _RS$  be the extended $S$-module of  an $R$-module $N_0$. There is the equality of groups
$$\Hr^0(k^G, M) = \Hr^0(G, \Aut _S(M))$$ and an isomorphism of pointed sets
$$ \Hr^1(k^G, M) \cong \Hr^1(G, \Aut _S(M)).$$ }

\goodbreak
\noindent{\sl  Proof  of Theorem 3.1.}
\smallskip
\noindent {\sl  $0$-level.} Let $\varphi$ be an element of   $\Hr^0(H, M)$. Then, by definition we have
$d^0\varphi = d^1\varphi$. This equality implies
$$(\id _M \ot \mu _H) (d^0\varphi \ot \id _H) \Delta _M = (\id _M \ot \mu _H) (d^1\varphi \ot \id _H) \Delta _M.$$
Let us compute the left-hand side on an element $m \in M$. We get
$$\eqalign{(\id _M \ot \mu _H) (d^0\varphi \ot \id _H) \Delta _M (m) &= \varphi(m_0)_0 \ot \varphi(m_0)_1\sigma_H(m_1)m_2
= \varphi(m_0)_0 \ot \varphi(m_0)_1\varepsilon_H(m_1)\cr & =
(\Delta _M \circ \varphi)(m_0\varepsilon_H(m_1)) = (\Delta _M \circ \varphi)(m).\cr}$$
The right-hand side applied to $m\in M$ is equal to
$$(\id _M \ot \mu _H) (d^1\varphi \ot \id _H) \Delta _M (m) = \varphi(m_0) \ot 1_Hm_1
= \varphi(m_0) \ot m_1  =
(\varphi \ot \id _H)\Delta _M (m).$$
Thus, one has $\Delta _M \circ \varphi = (\varphi \ot \id _H)\circ \Delta _M$, and therefore $f$ belongs to $\Dr^0(H, M)$.

\goodbreak
\smallskip
Conversely, let $f$ be an element of $\Dr^0(H, M)$. It satisfies the relation
$ (f \ot \id _H)\circ \Delta _M = \Delta _M \circ f $. 
Compose  each term of this equality  on the left with 
$(\id _M \ot \mu _H) \circ (\Delta _M \ot \sigma _H)$. The left-hand side becomes then exactly $d^0f$.
Apply the right-hand side on $m \in M$. Setting $m' = f(m)$, we get
$$m'_0\ot m'_1\sigma_H(m'_2) = m'_0 \ot \varepsilon_H(m'_1)1_H = m' \ot 1_H = f(m)\ot 1_H = d^1f(m).$$
Therefore $d^0f$ equals $d^1f$, hence $f$ belongs to $\Hr^0(H, M)$.

\medskip

\noindent  {\sl  $1$-level.} We begin to prove that $\tilde \kappa$ restricts to a bijection, still denoted by $\tilde \kappa$,
from  $ \Zr^1(H, M)$ to $\Cr^1(H, M)$. With the aim  to do that, we shall show that via  $\tilde \kappa$, for any $i = 1,2,3$,
Condition $\Zr\Cr_i$ of \pa 1.1 is equivalent to Condition $\Cr\Cr_i$ of \pa 2.1. We then prove that the bijection $\tilde \kappa$ induces a quotient map 
$\kappa:  \Hr^1(H, M) \lr \Dr^1(H, M)$ which is an isomorphism.
Adopt the following notations. For $\Phi \in \Zr^1(H, M)$ and $m \in M$, we  denote the tensor
$\Phi (m) \in M \ot H$ by $m_{[0]} \ot m_{[1]}$. Similarly, for $F \in \Cr^1(H, M)$ and $m \in M$, we set
$F (m) = m_{(0)} \ot m_{(1)}$.
\medskip
\goodbreak
\noindent {\sl  -- Equivalence of Condition $\Zr\Cr_1$ and Condition $\Cr\Cr_1$.} 
Fix $\Phi \in \Zr^1(H, M)$ and set $F = \tilde \kappa (\Phi) = \Phi \lodot \Delta _M$. So, for any $m\in M$, we have $F(m) = (m_0)_{[0]} \ot (m_0)_{[1]}m_1$.
Pick now $s \in S$. Condition $\Zr\Cr_1$ on $\Phi$ means $(ms)_{[0]} \ot (ms)_{[1]} = m_{[0]}s \ot m_{[1]}$. 
Let us compute  $F(ms)$:

$$\eqalign{F(ms) &= ((ms)_0)_{[0]} \ot ((ms)_0)_{[1]}(ms)_1\cr
&=  (m_0s_0)_{[0]} \ot (m_0s_0)_{[1]}m_1s_1 \cr
&=  (m_0)_{[0]} s_0\ot (m_0)_{[1]}m_1s_1\cr &= F(m) \Delta _S(s).\cr}$$
We use here the fact that $\Delta _M$ is twisted $S$-linear (second equality). Hence $F$ verifies Condition $\Cr\Cr_1$.
\smallskip

Conversely, fix $F \in \Cr^1(H, M)$.
Condition $\Cr\Cr_1$ on $F$ means 
$(ms)_{(0)} \ot (ms)_{(1)} = m_{(0)}s_0 \ot m_{(1)}s_1$, for any $s$  in $S$.
Set $\Phi = \tilde \kappa ^{-1}(F) = F \lodot \Delta _M'$, so
$\Phi(m) = (m_0)_{(0)} \ot (m_0)_{(1)}\sigma_H(m_1)$. Compute  $\Phi(ms)$:

$$\eqalign{\Phi(ms) &= ((ms)_0)_{(0)} \ot ((ms)_0)_{(1)}\sigma_H((ms)_1)\cr
&=  (m_0s_0)_{(0)} \ot (m_0s_0)_{(1)}\sigma_H(s_1)\sigma_H(m_1) \cr
&=  (m_0)_{(0)} s_0\ot (m_0)_{(1)}s_1\sigma_H(s_2)\sigma_H(m_1)\cr 
&= (m_0)_{(0)} s_0\ot (m_0)_{(1)}\varepsilon_H(s_1)\sigma_H(m_1)\cr
&= (m_0)_{(0)} s\ot (m_0)_{(1)}\sigma_H(m_1)\cr
&= \Phi(m)(s \ot 1).\cr}$$
Thus $\Phi$ verifies Condition  $\Zr\Cr_1$. 
\goodbreak

\medskip
\noindent {\sl  -- Equivalence of Condition $\Zr\Cr_2$ and Condition $\Cr\Cr_2$.} 
We still take $\Phi \in \Zr^1(H, M)$ and set $F = \tilde \kappa (\Phi) = \Phi \lodot \Delta _M$, so
$F(m) = (m_0)_{[0]} \ot (m_0)_{[1]}m_1$, for any $m\in M$.
Pick  $s \in S$. Condition
$\Zr\Cr_2$ on $\Phi$ is given by  the relation $m_{[0]}\varepsilon _H(m_{[1]}) = m$.
 Let us verify Condition $\Cr\Cr_2$ for $F$:
$$\eqalign{(\id _M \ot \varepsilon _H)F(m) &= (m_0)_{[0]}  \varepsilon _H((m_0)_{[1]})\varepsilon _H(m_1)\cr
&=  (m_0) \varepsilon _H(m_1)\cr
&=  (\id _M \ot \varepsilon _H)\Delta _M(m)\cr
&=  m.\cr}$$

Conversely, if $F$   verifies Condition $\Cr\Cr_2$, an easy computation shows that $\Phi =  F \lodot \Delta _M'$ fulfils Condition $\Zr\Cr_2$.

\medskip
\noindent {\sl  -- Equivalence of Condition $\Zr\Cr_3$ and Condition $\Cr\Cr_3$.} 
 We introduce the deformed differential map $\delta: \Wr_S^1(M) \lr \Wr_S^2(M)$
defined  on $\Phi \in \Wr_S^1(M)$ by the formula
$$\delta\Phi  = (\id _{M} \ot  T) \circ d^2\Phi$$ (recall that $T$ is the flip of $H\ot H$, see \pa 1.1).
%%Although we do not need this property, we observe  that $\delta$ is morphism
%%of monoids (straightforward calculation).
We prove now that Condition $\Cr\Cr_3$ on $F \in \Wr_S^1(M)$ may be translated into the equality
$$d^2F \lodot \delta F = d^1F. \eqno (4)$$ 
Indeed, as a consequence of the definitions of $\lodot$ and of $d^2$,
one gets $$d^2F \lodot \delta F = (\id _M \ot \mu _H ^{\otimes 2})(\id _M \ot \chi_2)(d^2F \ot \id _H^{\otimes 2})\delta F
= (\id _M \ot \mu _H ^{\otimes 2})(\id _M \ot \chi_2)(F\ot \eta_H \ot \id _H^{\otimes 2})\delta F.$$ 
Take $m \in M$ and observe that we have $ \delta F (m) = m_{(0)}\ot 1 \ot m_{(1)}$.  Let us compute
 $(d^2F \lodot \delta F) (m)$:
$$\eqalign{(d^2F \lodot \delta F)(m) 
& = (\id _M \ot \mu _H ^{\otimes 2})(\id _M \ot \chi_2)(F\ot \eta_H \ot \id _H^{\otimes 2})\delta F (m)\cr
& = (\id _M \ot \mu _H ^{\otimes 2})(\id _M \ot \chi_2)(F\ot \eta_H \ot \id _H^{\otimes 2})(m_{(0)}\ot 1 \ot m_{(1)})\cr
& = (\id _M \ot \mu _H ^{\otimes 2})(\id _M \ot \chi_2)(F(m_{(0)})\ot 1 \ot 1 \ot  m_{(1)})\cr
& = F(m_{(0)})\ot  m_{(1)}\cr
& = \bigl((F \ot \id _H) \circ F \bigr)(m).\cr}$$
Since $d^1F = (\id _M \ot \Delta _H)\circ F$, 
Condition $\Cr\Cr_3$ is equivalent to Equality (4).
\goodbreak
Let $\Phi$ be an element of $\Wr_S^1(M)$. Set $F = \tilde \kappa (\Phi) = \Phi \lodot \Delta _M$.
We write down a sequence of equivalent assertions which begins with Condition $\Zr\Cr_3$ on $\Phi$ and ends with an avatar of (4).

$$\eqalign{d^2\Phi \lodot d^0\Phi = d^1\Phi & \ \Longleftrightarrow  \
d^2(F \lodot \Delta _M') \lodot d^0(F \lodot \Delta _M') = d^1(F \lodot \Delta _M') \cr
& \ \Longleftrightarrow  \  
d^2F \lodot d^2\Delta _M' \lodot d^0F \lodot d^0\Delta _M' = d^1F \lodot d^1\Delta _M'\cr
& \ \Longleftrightarrow  \  
d^2F \lodot (d^2\Delta _M' \lodot d^0F \lodot d^0\Delta _M'\lodot d^1\Delta _M) = d^1F\cr}$$

It suffices now to prove  $d^2\Delta _M' \lodot d^0F \lodot d^0\Delta _M'\lodot d^1\Delta _M = \delta F$.
For any $m \in M$, one has the two equalities $d^0\Delta _M'(m) = m_0 \ot m_1\sigma_H(m_3)\ot \sigma_H(m_2)$
and $d^1\Delta _M(m) = m_0 \ot m_1 \ot m_2$.
Thus one gets $$(d^0\Delta _M'\lodot d^1\Delta _M)(m) = m_0 \ot m_1\sigma_H(m_3)m_4 \ot \sigma_H(m_2)m_5
= m_0 \ot m_1\ot \sigma_H(m_2)m_3 =  m_0 \ot m_1\ot 1. \eqno (5)$$
It remains to compute $(d^2\Delta _M' \lodot d^0F)(m)$.
Denote the tensor $d^0F(m_0) \in M\ot H$ by $x \ot y$, the summation being implicitly 
understood. Then $d^0F(m)$ is given by $ x_0 \ot x_1\sigma _H(m_1) \ot y$.
We also have $d^2\Delta _M'(m) = m_0\ot \sigma_H(m_1) \ot 1$. Therefore we get
 $$(d^2\Delta _M' \lodot d^0F)(m) = x_0 \ot \sigma_H(x_1)x_2\sigma_H(m_1)\ot 1y=
x \ot\sigma_H(m_1)\ot y. \eqno (6)$$
Combining (5) and (6), one obtains
$$\eqalign{((d^2\Delta _M' \lodot d^0F) \lodot (d^0\Delta _M'\lodot d^1\Delta _M)) (m) &= x \ot \sigma_H(m_1)m_2\ot y1\cr
&= x \ot \varepsilon_H(m_1)1\ot y\cr 
&= (\id _M \ot T)\bigl (x \ot y \ot \varepsilon_H(m_1)1\bigr)\cr
&= (\id _M \ot T) (F \ot \id _H)\bigl (m_0 \ot \varepsilon_H(m_1)1\bigr)\cr
&= (\id _M \ot T) (F(m) \ot 1)\cr
&= (\id _M \ot T) (d^2F)(m)\cr
&= (\delta F)(m).\cr} $$

\smallskip

\noindent {\sl  -- Factorization of $\tilde \kappa$.} We claim that the bijection $\tilde \kappa$ factorizes through an isomorphism from 
$\Hr^1(H, M)$ to $\Dr^1(H, M)$.
Indeed, take $\Phi$ and $\Phi '$ two cohomologous $1$-Hopf cocycles and $f \in \Aut _S(M)$ satisfying the equality
$d^1f^{-1} \lodot \Phi \lodot d^0f = \Phi '$.
Set $F = \tilde \kappa (\Phi)$ and  $F' = \tilde \kappa (\Phi') $.
One has then the equivalences
$$\eqalign{d^1f^{-1} \lodot \Phi \lodot d^0f = \Phi ' & \ \Longleftrightarrow  \   
d^1f^{-1} \lodot (F\lodot \Delta _M')\lodot d^0f = F'\lodot\Delta _M' \cr
& \ \Longleftrightarrow  \   
F\lodot \Delta _M' \lodot d^0f \lodot\Delta _M   = d^1f \lodot F' \cr
& \ \Longleftrightarrow  \   
F\lodot d^1f   = d^1f \lodot F' \cr
& \ \Longleftrightarrow  \   
F \circ f = (f \ot \id _H) \circ F' \cr}$$
The last equality means that $F$ and $F'$ are descent-cohomologous.
Observe that the third equivalence is a consequence of the equality
$d^0f = \Delta _M \lodot d^1f \lodot \Delta _M'$, which may be easily checked by the reader. \dm

\bigskip

\noindent {\sl  Post-scriptum.} The present work in its first preprint version led T. Brzezi\'nski to generalize the descent cohomology to the coring framework [2].
For any coring $C$ and any {\hbox{ $C$-comodule}} $M$, this author defines 
two descent cohomology sets $D^0(C, M)$  and $D^1(C, M)$, which coincide respectively with 
$D^0(H, M)$  and $D^1(H, M)$   (notations of \pa 2) when $C$ is 
the coring $S \ot H$.

\bigskip
\bigskip
\bigskip
\bigskip

{\ninerom

\noindent {R{\eightcmr EFERENCES}}

\medskip
\medskip

\item{[1]} A. B{\eightcmr LANCO} F{\eightcmr ERRO}, 
Hopf algebras and Galois descent, 
{\ninecmsl Publ. Sec. Mat. Universitat  Aut\`onoma Barcelona} {\ninecmbx 30} ({\oldstylen 1986}), no. 1, 65 -- 80.
\medskip

\item{[2]} T. B{\eightcmr RZEZI\'NSKI},
Descent cohomology and corings, 
Preprint arXiv: {\nineromtt math.RA/0601491} ({\oldstylen 2006}).
\medskip

\item{[3]} M. C{\eightcmr IPOLLA},  Discesa fedelmente piatta dei moduli,
{\ninecmsl Rendiconti del Circolo Matem\`atico di Palermo},
Serie II - tomo XXV ({\oldstylen 1976}).
\medskip

\item{[4]} P. D{\eightcmr ELIGNE}, Cat\'egories tannakiennes,
{\ninecmsl  The Grothendieck Festschrift}, Vol. II, 111 -- 195, 
Progr. Math., 87, 
Birkh\"auser Boston, Boston, MA ({\oldstylen 1990}). 
\medskip

\item{[5]} Y. D{\eightcmr OI},  M. T{\eightcmr AKEUCHI}, 
Hopf-Galois extensions of algebras, the Miyashita-Ulbrich action, and Azumaya algebras,
{\ninecmsl J. Algebra} {\ninecmbx 121} ({\oldstylen 1989}), no. 2, 488 -- 516.
\medskip

\item{[6]} M.A. K{\eightcmr NUS},  
{\ninecmsl  Quadratic and hermitian forms over rings},
 Grundlehren der mathematischen Wissenschaften 294, 
Springer-Verlag, Berlin -
Heidelberg - New York ({\oldstylen2 1991}).
\medskip

\item{[7]} H. F. K{\eightcmr REIMER}, M. T{\eightcmr AKEUCHI}, Hopf algebras and 
Galois extensions of an algebra,
{\ninecmsl Indiana Univ. Math. J.} %Vol. 
{\ninecmbx 30} ({\oldstylen 1981}), no. 5,
675  --  692.
\medskip

\item{[8]}  S. L{\eightcmr ANG}, J. T{\eightcmr ATE}, Principal homogeneous spaces over abelian varieties,
{\ninecmsl Amer. J. Maths.}   {\ninecmbx 80} ({\oldstylen 1958}), 
{\hbox{659  --  684.}}
\medskip

\item{[9]} L. {\eightcmr LE} B{\eightcmr RUYN}, M. {\eightcmr VAN DEN}  B{\eightcmr ERGH}, 
F. {\eightcmr VAN} O{\eightcmr YSTAEYEN},
{\ninecmsl Graded orders}, Birkh\"auser, Boston -- Basel ({\oldstylen 1988}).
\medskip

\item{[10]} P. N{\eightcmr USS},   Noncommutative descent and non-abelian cohomology, 
{\ninecmsl {\nineromi K}-Theory} {\ninecmbx 12} ({\oldstylen 1997}), no. 1, 23 -- 74.
\medskip

\item{[11]}  H.-J. S{\eightcmr CHNEIDER}, Principal homogeneous spaces
for arbitrary Hopf algebras,
{\ninecmsl Israel J. Math.} {\ninecmbx 72} ({\oldstylen 1990}), {\hbox {no. 1~--~2,}} 
167 -- 195.
\medskip

\item{[12]} J.-P. S{\eightcmr ERRE}, 
{\ninecmsl  Corps locaux}, Troisi\`eme \'edition corrig\'ee,
Hermann, Paris  ({\oldstylen 1968}).

\medskip

\item{[13]} J.-P. S{\eightcmr ERRE}, 
{\ninecmsl  Galois cohomology}, Springer-Verlag, Berlin --
Heidelberg  ({\oldstylen 1997}). Translated from
{\ninecmsl  Cohomologie galoisienne},
Lecture Notes in Mathematics 5, Springer-Verlag, Berlin --
Heidelberg -- New York ({\oldstylen 1973}).

\medskip

\item{[14]} M. E. S{\eightcmr WEEDLER},  
 Cohomology of algebras over Hopf algebras, {\ninecmsl Trans. Amer. Math. Soc.} {\ninecmbx 133}  ({\oldstylen 1968}), 205~--~239.

\vfill
\eject

\bye